\documentclass[12pt,a4paper]{article}

\usepackage{hyperref} 
\usepackage{amsmath,amsthm,amssymb,latexsym} 
\usepackage{amsfonts} 
\usepackage[american]{babel} 
\usepackage{mathrsfs} 
\usepackage[dvips]{graphicx} 
 
\newtheorem{teo}{Theorem}
\newtheorem{prop}[teo]{Proposition} 
\newtheorem{lemma}[teo]{Lemma}

\newenvironment{dimo}{{\it Proof.}}{$\square$\vspace{0.3cm}}
\newcommand{\R}{{\bf R}} 
 
\newcommand{\sfe}{{\bf S}^{n-1}}

\begin{document} 
\begin{title}
{From the Brunn-Minkowski inequality\\
to a class of Poincar\'e type inequalities}
\end{title} 
 
\author{Andrea Colesanti
\footnote{Dipartimento di Matematica ``U. Dini'', 
viale Morgagni 67/A, 50134 Firenze, Italy;\newline
colesant@math.unifi.it}} 
\date{} 
\maketitle

\begin{abstract} 
\noindent We present an argument which leads from the Brunn-Minkowski
inequality to a Poincar\'e type inequality on the boundary of a convex
body $K$ of class $C^2_+$ in ${\bf R}^n$. We prove that for every $\psi\in
C^1(\partial K)$  
\begin{equation*}
\int_{\partial K}\psi\,d{\cal H}^{n-1}=0\,\Rightarrow\,
\int_{\partial K}{\rm tr}(D\nu_K)\psi^2\,d{\cal H}^{n-1}
\le
\int_{\partial K}((D\nu_K)^{-1}\nabla\psi,\nabla\psi)\,d{\cal H}^{n-1}\,.
\end{equation*}
Here  ${\cal H}^{n-1}$ denotes the $(n-1)$-dimensional Hausdorff measure, 
$\nu_K$ is the Gauss map of $K$ and $D\nu_K$ is the differential of $\nu_K$
i.e. the Weingarten map. 

\end{abstract} 
 
\bigskip 

\noindent
{\em AMS 2000 Subject Classification:} 52A20; 26D10.

\section{Introduction} 

The Brunn-Minkowski inequality asserts that if $K_0$ and $K_1$ are
convex bodies (non-empty compact convex sets) in $\R^n$ and
$t\in[0,1]$ then  
\begin{equation}
\label{E0.1}
[V_n((1-t)K_0+tK_1)]^{1/n}\ge
(1-t)[V_n(K_0)]^{1/n}+t[V_n(K_1)]^{1/n}\,,
\end{equation}
where $V_n$ denotes the $n$-dimensional volume. 
This inequality is among the fundamental results in convex geometry,
i.e. the theory of convex bodies. More precisely, (\ref{E0.1}) can be
considered the core of the part of convex geometry concerning geometric
inequalities; it is the starting point towards many other similar inequalities
involving mixed volumes of convex bodies and it is related, as a special case,
to the Aleksandrov-Fenchel inequalities, another milestone in the theory of
convex bodies. We refer the reader to the book by Schneider \cite{Schneider},
and in particular to Chapter 6, for a detailed presentation of the
Brunn-Minkowski inequality in the context of convex geometry.

However, as it is underlined in the survey paper \cite{Gardner} by Gardner,
the r\^ole of (\ref{E0.1}) goes beyond the boundaries of the theory of convex
bodies. The main evidence of this fact are its connections 
with several important inequalities in
analysis. The Brunn-Minkowski inequality provides a simple proof of the
isoperimetric inequality in the case of smooth domains. Moreover, it admits an
equivalent functional formulation, the Pr\'ekopa-Leindler inequality, which is
related to Young's convolution inequality (or better, to its reverse
form). All these links, and many others, are very well described in \cite{Gardner}. 
Recently, Bobkov and Ledoux in \cite{Bobkov-Ledoux2} gave a proof based on the
Brunn-Minkowski inequality of the Sobolev and Gagliardo-Nirenberg inequalities
with optimal constant. The approach used in \cite{Bobkov-Ledoux2} is analogous to the one
previously used by the same authors in \cite{Bobkov-Ledoux1}, where they 
prove the logarithmic Sobolev
inequality and an inequality of Poincar\'e type due Brascamp and Lieb (see 
\cite{Brascamp-Lieb}); we will come back to the latter inequality at the end
of this introduction.

In this paper we present an argument which leads from the Brunn-Minkowski
inequality to a Poincar\'e type inequality on the boundary of smooth convex
bodies with positive Gauss curvature. Let us state our main result.
A convex body $K$ in ${\bf R}^n$ is said to be of class $C^2_+$ if its
boundary $\partial K$ is of class $C^2$ and the Gauss curvature is strictly
positive at each point of $\partial K$.  

\begin{teo} 
\label{T1} Let $K\subset\R^n$ be a convex body of class $C^2_+$ and let $\nu_K$
be the Gauss map of $K$. For every $\psi\in C^1(\partial K)$, if 
\begin{equation}
\label{E0.2}
\int_{\partial K}\psi\,d{\cal H}^{n-1}=0\,,
\end{equation}
then
\begin{equation}
\label{E0.3}
\int_{\partial K}{\rm tr}(D\nu_K)\psi^2\,d{\cal H}^{n-1}
\le
\int_{\partial K}((D\nu_K)^{-1}\nabla\psi,\nabla\psi)\,d{\cal H}^{n-1}\,.
\end{equation}
Here ${\cal H}^{n-1}$ denotes the $(n-1)$-dimensional Hausdorff measure and
$D\nu_K$ is the differential of the Gauss map, i.e. the Weingarten map.
\end{teo}

As a first comment we observe that inequality (\ref{E0.3}) is sharp; indeed we will see in \S
\ref{Paragrafo3} that if $\psi(x)=(\nu_K(x),u_0)$, $x\in\partial K$, where $u_0$ is a fixed
vector in $\R^n$, then (\ref{E0.2}) is fulfilled and (\ref{E0.3}) becomes an
equality. 

If we choose $K$ to be the unit ball, then $\partial K=\sfe$ and $\nu_K$ is
the identity map on $\sfe$. In this case from Theorem \ref{T1} we get that
\begin{equation*}
(n-1)
\int_{\sfe}\psi^2\,d{\cal H}^{n-1}
\le
\int_{\sfe}|\nabla\psi|^2\,d{\cal H}^{n-1}\,,
\end{equation*}
for every $\psi\in C^1(\sfe)$ such that
\begin{equation*}
\int_{\sfe}\psi\,d{\cal H}^{n-1}=0\,.
\end{equation*}
This is the classic Poincar\'e inequality on $\sfe$ with the sharp constant.

Let $K$ be as in Theorem \ref{T1} and let $\alpha>0$ be such that all the
principal curvatures of $\partial K$ are greater than or equal to $\alpha$. We recall that
the eigenvalues of $D\nu_K$ are just the principal curvatures of $\partial K$, hence
$$ 
{\rm tr}(D\nu_K)\ge(n-1)\alpha\quad{\rm and}\quad
((D\nu_K)^{-1}\nabla\psi,\nabla\psi)\le\frac{1}{\alpha}|\nabla\psi|^2\quad{\rm
  on}\;\partial K\,.
$$
Then from Theorem \ref{T1} we deduce that
\begin{equation}
\label{E0.4}
(n-1)\alpha^2
\int_{\partial K}\psi^2\,d{\cal H}^{n-1}
\le
\int_{\partial K}|\nabla\psi|^2\,d{\cal H}^{n-1}\,,
\end{equation}
for all $\psi\in C^1(\partial K)$ such that (\ref{E0.2}) holds. Inequality
(\ref{E0.4}) provides a lower bound for the first eigenvalue $\lambda(\partial
K)$ of the Laplace operator on $\partial K$:
\begin{equation*}
  \lambda(\partial K)\ge (n-1)\alpha^2\,.
\end{equation*}
This is a special case of a theorem due to Lichnerowicz which gives a similar lower
bound in the case of Riemannian manifolds verifying a suitable assumptions on the Ricci
tensor (see \cite{Chavel}, Chapter 3). 

We also would like to point out the analogy between Theorem \ref{T1} and an
inequality due to Brascamp and Lieb (see Theorem 4.1 in \cite{Brascamp-Lieb}
and \cite{Bobkov-Ledoux2} for a proof based on the Pr\'ekopa-Leindler inequality). 
Let $f=e^{-u}$ where $u\in C^2({\bf R}^n)$, $D^2u>0$ in ${\bf R}^n$ and
$\lim_{|x|\to\infty}u(x)=\infty$. Define the measure $\mu$ in ${\bf R}^n$
through the equality
$$
d\mu=f\,dx\,.
$$
Then for every $\psi\in C^1({\bf R}^n)$ having compact support,
\begin{equation}
\label{E0.5}
\int_{{\bf R}^n}\psi\,d\mu=0\;
\Rightarrow\;\int_{{\bf R}^n}\psi^2\,d\mu\le
\int_{{\bf R}^n}((D^2u)^{-1}\nabla\psi,\nabla\psi)\,d\mu\,.
\end{equation}
Like (\ref{E0.3}), (\ref{E0.5}) can be used to prove inequalities of Poincar\'e
type. Choosing  
$u(x)=\frac{\Vert x\Vert^2}{2}$, $x\in\R^n$, (\ref{E0.5}) becomes the Poincar\'e inequality
for the Gaussian measure. More generally, if $D^2u\ge c\,{\it Id}$, where ${\it
  Id}$ is the identity matrix and $c>0$, then (\ref{E0.5}) provides a
Poincar\'e inequality for the measure $\mu$ (see also \cite{Bobkov} for  
Poincar\'e inequalities fulfilled by general log-concave measures).

The idea of the proof of Theorem \ref{T1} can be heuristically described in
the following way. The Brunn-Minkowski inequality asserts that the
volume $V_n$ raised to the power $1/n$ is concave in the set
of $n$-dimensional convex bodies. In the subset of convex bodies of class $C^2_+$  
one can actually compute the {\em second variation} of $V_n^{1/n}$, and this has to
be negative semi-definite; this property is equivalent to Theorem
\ref{T1}. This approach was suggested to the author by the paper
\cite{Jerison2} of Jerison (see also \cite{Jerison1}), where the second
variation of the Newtonian capacity of convex bodies is determined. 

\section{Preliminaries}
In this section we recall some basic facts about convex bodies. The general
reference for the results reported here is \cite{Schneider} and in
particular \S 2.5 for the properties of convex bodies
of class $C^2_+$.

Let $K\subset\R^n$ be a convex body; we denote by $h_K$ the {\em support
  function} of $K$:
$$
h_K\,:\,\sfe\longrightarrow\R\,,\quad
h_K(u)=\sup\{(x,u)\,|\,x\in K\}\,.
$$
A convex body $K$ is said to be class $C^2_+$ if $\partial K\in C^2$ and the
Gauss curvature is strictly positive at each point of $\partial K$. In the
rest of this section we will assume that $K$ is of class $C^2_+$. We denote by 
$\nu_K\,:\,\partial K\rightarrow\sfe$ the Gauss map of $K$, i.e. for
$x\in\partial K$, $\nu_K(x)$ is the outer unit normal to $\partial K$ at
$x$. The map $\nu_K$ is
differentiable on $\partial K$ and its differential $D\nu_K$ is the Weingarten
map. The present assumptions imply that $\nu_K$ is invertible and the
inverse map is differentiable on $\sfe$.

For a function $f\in
C^2(\sfe)$ we denote by $f_i$ and $f_{ij}$, $i,j=1,\dots,n-1$, the first and
second covariant derivatives of $f$ with respect to an orthonormal frame
$\{e_1,\dots,e_{n-1}\}$ on $\sfe$. 

As $K$ is of class $C^2_+$, $h_K\in
C^2(\sfe)$. The following matrix will play an important r\^ole throughout this
paper: 
\begin{equation*}
\label{E1.1}
((h_K)_{ij}+h_K\delta_{ij})\,,
\end{equation*}
where $\delta_{ij}$, $i,j=1,\dots,n-1$, are the standard Kronecker symbols.
Let $u\in\sfe$, the linear map induced by $((h_K)_{ij}(u)+h_K(u)\delta_{ij})$
on the tangent space of $\sfe$ at $u$ is the {\em reverse Weingarten map} of
$K$ at $u$. In other words
\begin{equation}
  \label{E1.2}
((h_K)_{ij}+h_K\delta_{ij})=D(\nu_K^{-1})\quad{\rm on}\;\sfe\,.
\end{equation}
Moreover we have the matrix inequality
\begin{equation}
  \label{E1.3}
((h_K)_{ij}+h_K\delta_{ij})>0\quad{\rm on}\;\sfe\,;  
\end{equation}
indeed for $u\in\sfe$ the eigenvalues of $((h_K)_{ij}(u)+h_K(u)\delta_{ij})$ are the
principal radii of curvature of $K$ at the point $\nu^{-1}(u)$.  Conversely,
let $h\in C^2(\sfe)$ be such that $(h_{ij}+h\delta_{ij})>0$ on $\sfe$; then
there exists a unique convex body $K$ of class $C^2_+$ such that $h=h_K$. 
 
The following formula expresses the volume of $K$ in
terms of its support function:
\begin{equation}
  \label{E1.4}
V_n(K)=\frac{1}{n}\int_{\sfe}h_K\,{\rm det}((h_K)_{ij}+h_K\delta_{ij})\,d{\cal H}^{n-1}\,.  
\end{equation}

We introduce a set of functions defined on $\sfe$:
$$
{\cal C}=\{h\in C^2(\sfe)\,:\,(h_{ij}+h\delta_{ij})>0\;{\rm on}\;\sfe\}\,.
$$
By the previous remarks $\cal C$ is formed by the support
functions of convex bodies of class $C^2_+$. Next we define the functional
\begin{equation}
\label{E1.4b}
F\,:\,{\cal C}\longrightarrow\R_+\,,\quad
F(
h)=\frac{1}{n}\int_{\sfe}h\,{\rm
  det}(h_{ij}+h\delta_{ij})\,d{\cal H}^{n-1}\,,
\end{equation}
i.e. $F$ is the volume functional. 

One of the basic properties of support
functions is that if $K_0$ and $K_1$ are convex bodies and $\alpha,\beta\ge0$,
then $h_{\alpha K_0+\beta K_1}=\alpha h_{K_0}+\beta h_{K_1}$. The following proposition
is a direct consequence of this fact and of the Brunn-Minkowski inequality
(\ref{E0.1}).  

\begin{prop} $F^{1/n}$ is concave in $\cal C$, i.e. 
\label{P1}
\begin{equation}
\label{E1.5}
F^{1/n}((1-t)h_0+th_1)\ge
(1-t)F^{1/n}(h_0)+tF^{1/n}(h_1)\,,\quad
\forall\,h_0,h_1\in{\cal C}\,,\;\forall\,t\in[0,1]\,.
\end{equation}
\end{prop}


In the next section we will use the notion of {\em cofactor matrix}, that we
briefly recall. Let $A=(a_{ij})$ be a $k\times k$ matrix; for every $i,j=1,\dots,k$,
define
\begin{equation}
\label{E1.6}
\gamma_{ij}=\frac{\partial{\rm det}(A)}{\partial a_{ij}}\,.
\end{equation}
The matrix $(\gamma _{ij})$ is called the cofactor matrix of $A$. If $A$ is
invertible then $(\gamma_{ij})={\rm det}(A)\, A^{-1}$. As a simple consequence of
the elementary properties of the determinant we note the following formula:
\begin{equation}
  \label{E1.7}
\sum_{i,j=1}^k\gamma_{ij}\,a_{ij}=k\det(A)\,.  
\end{equation}

For an arbitrary $h\in{\cal C}$ we denote by $(c_{ij})$ the cofactor matrix of  
$(h_{ij}+h\delta_{ij})$. The proof of the next lemma can be found in
\cite{Cheng-Yau} (see page 504). 

\begin{lemma} For every $j=1,\dots,n-1$
\label{L1}
$$
{\rm div}_i(c_{ij})=\sum_{i=1}^{n-1}(c_{ij})_i=0\,.
$$  
\end{lemma}

\section{The proof of Theorem \ref{T1}}
\label{Paragrafo3}

We start with the following theorem.

\begin{teo}
  \label{T2}
Let $h\in{\cal C}$, $\phi\in C^1(\sfe)$ and denote by $(c_{ij})$ the cofactor matrix of  
$(h_{ij}+h\delta_{ij})$. If
\begin{equation}
\label{E2.1}
\int_{\sfe}\phi\det(h_{ij}+h\delta_{ij})d{\cal H}^{n-1}=0
\end{equation}
then
\begin{equation}
\label{E2.2}
\int_{\sfe}{\rm tr}(c_{ij})\phi^2d{\cal H}^{n-1}\le
\int_{\sfe}\sum_{i,i=1}^{n-1}c_{ij}\phi_i\phi_jd{\cal H}^{n-1}\,.
\end{equation}
\end{teo}

As we will see, Theorem \ref{T1} and Theorem \ref{T2} can
be obtained one from each other through the change of variable given by the Gauss
map of $K$, where $K$ is such that $h=h_K$, and its inverse. 

The proof of Theorem \ref{T2} is based on a simple idea. Assume for a moment 
that $\phi\in C^\infty(\sfe)$; for $s\in\R$ with $|s|$ sufficiently small, $h+s\phi\in{\cal
  C}$. By Proposition \ref{P1} the function of one variable $g(s)=[F(h+s\phi)]^{1/n}$
is concave. Inequality (\ref{E2.2}), under the condition (\ref{E2.1}), will
follow simply imposing $g''(0)\le 0$.

\begin{prop}\label{P2}
Let $h\in{\cal C}$, $\phi\in C^\infty(\sfe)$ and $\varepsilon>0$ be such that
$h+s\phi\in{\cal C}$ for every $s\in(-\varepsilon,\varepsilon)$. Set 
$h_s=h+s\phi$ and $f(s)=F(h_s)$, where $F$ is defined by (\ref{E1.4b}). Then 
\begin{equation}
  \label{E2.3}
f'(s)=\int_{\sfe}\phi\det((h_s)_{ij}+h_s\delta_{ij})d{\cal H}^{n-1}\,,\quad
\forall\,s\in(-\varepsilon,\varepsilon)\,.
\end{equation}
\end{prop}

\noindent
\begin{dimo} For every $u\in\sfe$
\begin{eqnarray*}
&&\dfrac{d}{ds}[h_s(u)\det((h_s)_{ij}(u)+h_s(u)\delta_{ij})]=\\
&&=\phi(u)\det((h_s)_{ij}(u)+h_s(u)\delta_{ij})+
h_s(u)\sum_{i,j=1}^{n-1}c^s_{ij}(u)(\phi_{ij}(u)+\phi(u)\delta_{ij})\,,
\end{eqnarray*}
where $(c^s_{ij})$ denotes the cofactor matrix of
$((h_s)_{ij}+h_s\delta_{ij})$. Differentiating under the integral sign we obtain
\begin{equation}
\label{E2.3bis}
f'(s)=\frac{1}{n}\int_{\sfe}
\left[\phi\det((h_s)_{ij}+h_s\delta_{ij})+h_s\sum_{i,j=1}^{n-1}c^s_{ij}(\phi_{ij}+\phi\delta_{ij})\right]
d{\cal H}^{n-1}\,.
\end{equation}
Now
\begin{equation}
\label{E2.4}
\int_{\sfe}h_s\sum_{i,j=1}^{n-1}c^s_{ij}\phi_{ij}d{\cal H}^{n-1}=
\int_{\sfe}\phi\sum_{i,j=1}^{n-1}c^s_{ij}(h_s)_{ij}d{\cal H}^{n-1}\,,
\end{equation}
where we have integrated by parts twice and used Lemma \ref{L1}. On the other
hand by (\ref{E1.7}) we have
\begin{equation}
\label{E2.5}
\sum_{i,j=1}^{n-1}c^s_{ij}((h_s)_{ij}+h_s\delta_{ij})=(n-1)\det((h_s)_{ij}+h_s\delta_{ij})\,. 
\end{equation}
Inserting (\ref{E2.4}) and (\ref{E2.5}) into (\ref{E2.3bis}) completes the proof.
\end{dimo}

The next result is a straightforward consequence of Proposition \ref{P2} and
the definition of cofactor matrix.

\begin{prop}\label{P3} 
In the assumptions and notations of Proposition \ref{P2}
\begin{equation}
\label{E2.6}
f''(0)=\int_{\sfe}\phi\sum_{i,j=1}^{n-1}c_{ij}(\phi_{ij}+\phi\delta_{ij})\,d{\cal H}^{n-1}\,,
\end{equation}
where $(c_{ij})$ is the cofactor matrix of $(h_{ij}+h\delta_{ij})$.
\end{prop}

\bigskip

\noindent
{\it Proof of Theorem \ref{T2}.} Assume first that $\phi\in
C^\infty(\sfe)$. Let $\varepsilon>0$ be such that $h+s\phi\in{\cal C}$ for  
every $s\in(-\varepsilon,\varepsilon)$. Set $f(s)=F(h+s\phi)$ and
$g(s)=[f(s)]^{1/n}$ for $s\in(-\varepsilon,\varepsilon)$. By Proposition
\ref{P1} $g$ is concave, so that
\begin{equation*}
g''(0)=\frac{1}{n}\left(\frac{1}{n}-1\right)[f(0)]^{\frac{1}{n}-2}[f'(0)]^2+
\frac{1}{n}[f(0)]^{\frac{1}{n}-1}f''(0)\le0\,.  
\end{equation*}
In particular, if (\ref{E2.1}) holds then by Proposition \ref{P2} $f'(0)=0$ so
that $f''(0)\le 0$ (note that $f(0)=F(h)>0$). By Proposition \ref{P3} we
obtain
$$
\int_{\sfe}{\rm tr}(c_{ij})\phi^2\,d{\cal H}^{n-1}\le
-\int_{\sfe}\phi\sum_{i,j=1}^{n-1}c_{ij}\phi_{ij}\,d{\cal H}^{n-1}\,.
$$
We get inequality (\ref{E2.2}) integrating by parts the integral in the right
hand-side (and using Lemma \ref{L1}). The general case $\phi\in
C^1(\sfe)$ follows by a standard approximation argument. 
$\square$
\vspace{0.3cm}

\noindent
{\bf Remark.} Let $h\in{\cal C}$ and $\phi(u)=(u,u_0)$, $u\in\sfe$, where $u_0$ is a
fixed vector in $\R^n$. Let $K$ be the convex body such that $h_K=h$. For
every $s\in\R$ we have $h+s\phi=h_{K+su_0}$, i.e. $h+s\phi$ is the support
function of a translate of $K$. This implies in particular that
$$
F(h+s\phi)=V_n(K+su_0)=V_n(K)\quad\forall\, s\in\R\,.
$$ 
Then, in the notations of the proof of Theorem \ref{T2} we have
$f'(0)=0$ and $f''(0)=0$. These equalities tell us that (\ref{E2.1}) is
fulfilled and that (\ref{E2.2}) is an equality with this choice of the
function $\phi$.
\bigskip

\noindent
{\it Proof of Theorem \ref{T1}.} For simplicity
we set $h=h_K$ and $\nu=\nu_K$. As before, denote by $(c_{ij})$ the cofactor matrix of
$(h_{ij}+h\delta_{ij})$. Note that for every $u\in\sfe$
\begin{equation}
\label{E2.6bis}
(c_{ij}(u))=\det(h_{ij}(u)+h\delta_{ij}(u))\,(h_{ij}(u)+h\delta_{ij}(u))^{-1}\,. 
\end{equation}

Let $\psi\in C^1(\partial K)$ and set 
$\phi(u)=\psi(\nu^{-1}(u))$, $u\in\sfe$; then $\phi\in C^1(\sfe)$.
Performing the change of variable $u=\nu(x)$ and using (\ref{E1.2}) we obtain
\begin{equation}
\label{E2.7}
\int_{\sfe}\phi(u)\det(h_{ij}(u)+h\delta_{ij}(u))d{\cal H}^{n-1}(u)=
\int_{\partial K}\psi(x)\,d{\cal H}^{n-1}(x)\,.
\end{equation}
Analogously, using (\ref{E1.2}) and (\ref{E2.6bis}) we have
\begin{equation}
  \label{E2.8}
\int_{\sfe}{\rm tr}(c_{ij}(u))\phi^2(u)d{\cal H}^{n-1}(u)=
\int_{\partial K}{\rm tr}(D\nu_K(x))\psi^2(x)\,d{\cal H}^{n-1}(x)\,,  
\end{equation}
and, for every $u\in\sfe$,
\begin{eqnarray*}
\sum_{i,i=1}^{n-1}c_{ij}(u)\phi_i(u)\phi_j(u)
&=&\det(h_{ij}(u)+h\delta_{ij}(u))\left(D\nu^{-1}(u)\nabla\phi(u),\nabla\phi(u)\right)\\
&=&\det(h_{ij}(u)+h\delta_{ij}(u))\left(\nabla\psi(\nu^{-1}(u)),D\nu^{-1}(u)\nabla\psi(\nu^{-1}(u))\right)\,.
\end{eqnarray*}
Then
\begin{equation}
\label{E2.9}
\int_{\sfe}\sum_{i,i=1}^{n-1}c_{ij}(u)\phi_i(u)\phi_j(u)d{\cal H}^{n-1}(u)=
\int_{\partial K}((D\nu(x))^{-1}\nabla\psi(x),\nabla\psi(x))\,d{\cal H}^{n-1}(x)\,.
\end{equation}
Theorem \ref{T1} follows from (\ref{E2.7}), (\ref{E2.8}), (\ref{E2.9}) and Theorem \ref{T2}.
$\square$
\vspace{0.3cm}

By the remark following the proof of Theorem \ref{T2}, if
$\psi(x)=(\nu_K(x),u_0)$, $x\in\partial K$, where $u_0$ is any fixed vector in
$\R^n$, then (\ref{E0.2}) is fulfilled and (\ref{E0.3}) becomes an equality.

\end{document}